\documentclass[12pt]{article}
\usepackage{graphicx}
\usepackage{epstopdf} 
\input rlepsf
\newtheorem{thm}{Theorem}[section]
\newtheorem{prop}[thm]{Proposition}

\newtheorem{lemma}[thm]{Lemma}
\newtheorem{cor}[thm]{Corollary}

\newtheorem{exer}[thm]{Exercise}
\usepackage{amsfonts} 
\usepackage{amsmath} 
\usepackage{amscd} 
\usepackage{amssymb}

\def\C{{\mathbb C}}
\def\a{{\alpha}}
\def\b{{\beta}}
\def\c{{\gamma}}
\def\s{{\sigma}}

\def\R{{\mathbb R}}
\def\T{{\mathcal T}}
\def\H{{\mathbb H}}
\def\Z{{\mathbb Z}}

\begin{document}

\begin{center}

{\Large Tutorial on the braid groups}

\end{center}

\vfill

Originally presented to the summer school and workshop
at the Institute for Mathematical Sciences, National University of Singapore, June  2007

\vfill

\begin{center}
{

Dale Rolfsen 

Department of Mathematics

University of British Columbia

Vancouver, BC, Canada

rolfsen@math.ubc.ca
}

\end{center}

\eject

\section*{Introduction.}  
The braid groups $B_n$ were introduced by E. Artin eighty years ago
 \cite{A1}, although their significance to mathematics was possibly realized a century earlier by Gauss, as evidenced by sketches of braids in his notebooks, and later in the nineteenth century by Hurwitz.  The braid groups provide a very attractive blending of geometry and algebra, and have applications
in a wide variety of areas of mathematics, physics, and recently in polymer chemistry and molecular biology.

Despite their ripe old age, and an enormous amount of attention paid
to them by mathematicians and physicists, these groups still provide us with big
surprises.  The subject was given a major boost when V. Jones discovered, in the mid 1980's,
new representations of the braid groups, which led to the famous Jones polynomial of knot theory.  In this minicourse, I will touch on that discovery, and more recent ones, regarding this
wonderful family of groups, and raise a few questions about them which are still
open.  Along the way, I will also discuss some properties of these groups which deserve to
be better known.  Among my personal favorites is the orderability of braid groups, which implies that they have some special algebraic properties.  That will be the subject of my last lecture.

Due to time restrictions, I will not be able to discuss many interesting aspects of braid theory, such as the conjugacy problem as solved by Garside, dynamics of braids, applications to cryptography and connections with homotopy groups of spheres.  Some of these issues will be considered in lectures by others in this summer school.  Also because of time limitation, many proofs will be omitted in my presentation.  Some proofs will be left as exercises for the students, often with hints or references.

The emphasis of my presentation is that there are many ways of looking at braids, and each point of view provides new information about these groups.  Though much of the motivation for studying braids comes from knot theory, my emphasis will be on the algebraic aspects of braid groups.  

\subsection*{Prerequisites}  Students will be assumed to have background in basic topology and group theory.  No other specialized expertise is necessary.

\eject

\subsection*{Suggested reading}

$\bullet$ Joan Birman and Tara Brendle, Braids, A Survey, in {\it Handbook of Knot Theory} (ed. W. Menasco and M. Thistlethwaite), Elsevier, 2005.

Available electronically at:

http://www.math.columbia.edu/~jb/Handbook-21.pdf

\bigskip

$\bullet$ Joan Birman, {\it Braids, links and mapping class groups}, Annals of Mathematics Studies, Princeton University Press, 1974.

This classic is still a useful reference.

\bigskip

$\bullet$ Gerhard Burde and Heiner Zieschang, {\it Knots}, De Gruyter Studies in Mathematics 5, second edition, 2003.

This contains a very readable chapter on braid groups.

\bigskip

$\bullet$ Patrick Dehornoy, Ivan Dynnikov, Dale Rolfsen, and Bert Wiest, {\it Why are braids orderable?}
Soc. Math. France series Panoramas et Syntheses 14(2002).  

Available electronically at:

http://perso.univ-rennes1.fr/bertold.wiest/bouquin.pdf

\bigskip

$\bullet$ Vaughan F. R. Jones, {\it Subfactors and Knots}, CBMS Regional Conference series in Mathematics, Amer. Math. Soc., 1991.

This gem focusses on the connection between operator algebras and knot theory, via the braid groups.

\bigskip

$\bullet$ Vassily Manturov, {\it Knot Theory}, Chapman \& Hall, 2004.

The most up-to date text.  Contains details, for example, of the faithfulness Bigelow-Krammer-Laurence representation.

$\bullet$ Kunio Murasugi and Bohdan Kurpita, {\it A study of braids},  Mathematics and its applications,
Springer, 1999.

\bigskip

Also see research papers and other books in the bibliography at the
end of the notes.

\eject

\section{Braids as strings, dances or symbols:}

One of the interesting things about the braid groups $B_n$ is they can be defined
in so many ways, each providing a unique insight.  I will briefly describe some
of these definitions, punctuated by facts about $B_n$ which are revealed from
the various points of view. You can check the excellent references \cite{A1}, \cite{A2},
\cite{Bir1}, \cite{BZ}, \cite{F}, \cite{K} for more information.

{\bf Definition 1: Braids as strings in 3-D.}  This is the usual and visually
appealing picture.  An $n$-braid is a collection of $n$ strings in
$(x,y,t)$-space, which are disjoint and monotone in (say) the $t$ direction.

We require that the endpoints of the strings are at fixed points, say the
points $(k, 0, 0)$ and $(k, 0, 1)$, $k = 1, \dots , n$.  We also regard two braids to be 
equivalent (informally, we will say they are equal) if one can be deformed to the 
other through the family of braids, with endpoints fixed throughout
the deformation.

Although most authors draw braids vertically, I prefer to view them
horizontally with the $t$-axis running from left to right.  My reason
is that we do algebra (as in writing) from left to right, and the
multiplication of braids is accomplished by concatenation (which is
made more explicit in the second definition below), which we take to
be in the same order as the product.  In the vertical point of view,
some authors view the product $\a\b$ of two braids $\a$ and $\b$ with
$\a$ above $\b$, while others would put $\b$ on top of $\a$.  The
horizontal convention eliminates this ambiguity.

There is, of course, a strong connection between braids and knots.  A braid
$\b$ defines a knot or link $\hat{\b}$, its closure, by connecting the 
endpoints in a standard way, without introducing further interaction 
between the strings.

\begin{center}%
\includegraphics[totalheight=6cm]{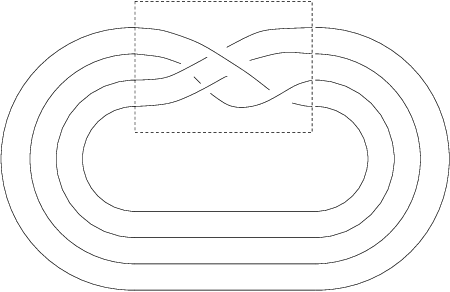}%
\begin{center}%
Figure 1:  The closure of a braid.
\end{center}
\end{center}

Equivalent braids give rise to equivalent links, but
different braids may give rise to the same knot or link.  We will
discuss how to deal with this ambiguity later.  $B_n$ will denote the
group formed by $n$-braids with concatenation as the product, which we will soon verify has the properties of a group, namely an associative multiplication, with identity element and inverses.

{\bf Definition 2: Braids as particle dances.}  If we take $t$ as a
time variable, then a braid can be considered to be the time history
of particles moving in the $(x,y)$ plane, or if you prefer, the
complex plane, with the usual notation  $x + y\sqrt{-1}$.  

  This gives the view of a braid as a dance of noncolliding particles
in $\C$, beginning and ending at the integer
points $\{ 1, \dots, n \}$. We think of the particles moving
in trajectories
$$\b(t) = (\b_1(t), \dots \b_n(t)),\quad \b_i(t) \in \C,$$ where $t$
runs from $0$ to $1$, and  $\b_j(t) \ne \b_k(t)$ when $j \ne k$.
  A {\it braid} is then
such a time history, or dance, of noncolliding particles in the plane
which end at the spots they
began, but possibly permuted. Equivalence of braids roughly reflects the 
fact that choreography does not specify precise positions of the
dancers, but rather their relative positions.  However, notice that if particle $j$ and $j+1$ interchange places, rotating clockwise, and then do the same move, but counterclockwise, then their dance is equivalent to just standing still!

The product of braids can be regarded as one dance following the other, both at double
speed.  Formally, if $\a$ and $\b$ are braids, we define their product $\a\b$
(in the notation of this definition) to be the braid which is $\a(2t)$
for $0 \le t \le 1/2$ and $\b(2t-1)$ for $1/2 \le t \le 1$.

\begin{exer} Verify that the product is associative, up to equivalence.  The identity braid is to stand still, and each dance has an inverse; doing the dance in reverse.  More formally, if $\b$ is a braid, as defined in this section, define $\bar{\b}(t) = \b(1-t)$.
Write a formula for the product $\c = \b\bar{\b}$ and verify that there is a continuous deformation 
from $\c$ to the identity braid.  [Hint: let $s$ be the deformation parameter, and consider the dance
which is to perform $\b$ up till time $s/2$, then stand still until time $1 - s/2$, then do the dance $\bar{\b}$ in the remaining time.]
\end{exer}

A braid $\b$ defines a permutation $i \to \b_i(1)$ which is a
well-defined element of the permutation group $\Sigma_n$.  This is a
homomorphism  $B_n \to \Sigma_n$
with kernel, by definition, the subgroup $P_n < B_n$ of {\it pure} braids. 
$P_n$ is sometimes called the {\it colored} braid group, as the particles can be
regarded as having identities, or colors. 
$P_n$ is of course normal in $B_n$, of index $n!$, and there is an exact sequence
$$1 \to P_n \to B_n \to \Sigma_n \to 1.$$

\begin{exer}
Show that any braid is equivalent to a piecewise-linear braid.  Moreover, one may assume that under the projection $p: \R^3 \to \R^2$ given by $p(x,y,t) = (y,t)$, there are only a finite number of 
singularities of $p$ restricted to the (image of the) given braid.  One may assume these are all
double points, where a pair of strings intersect transversely, and at different $t$-values.
\end{exer}

\begin{center}%
\includegraphics[totalheight=6cm]{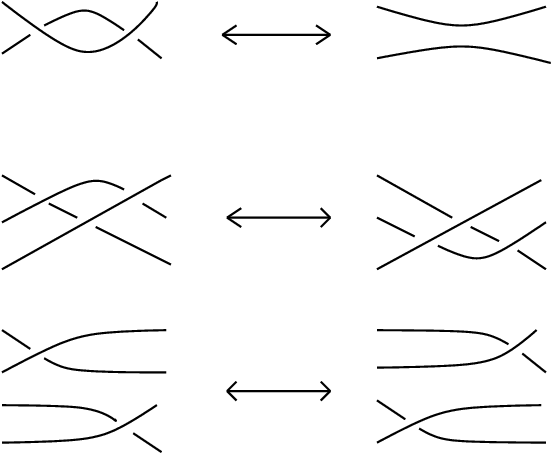}%
\begin{center}%
Figure 2:  Equivalent braids.
\end{center}
\end{center}

\begin{exer}  Viewing the projection in the $(y,t)$ - plane, with the $t$-axis as horizontal, for any fixed value of $t$ at which a crossing does not occur,
label the strings $1, 2, \dots ,n$ counting from the bottom to the top in the $y$ direction.  Let
$\sigma_i$ for $i = 1, \dots , n-1$ denote the braid with all strings horizontal, except that the string $i+1$ crosses over the $i$ string, resulting in the permutation $i \leftrightarrow i+1$ in the labels of the strings.  (Think of it as a right-hand screw motion of these two strings.  Argue that any braid is equivalent to a product of these ``generators'' of the braid group.  [In Figure 1, for example, the braid could be written either as $\s_2^{-1}\s_1^{-1}\s_3\s_2\s_1^{-1}$ or
$\s_2^{-1}\s_3\s_1^{-1}\s_2\s_1^{-1}$ (or indeed many other ways).]

Finally, note that the inverse of $\sigma_i$ is the braid with all strings horizontal, except that the string $i$ crosses over the $i+1$ string.  The first sketch in Figure 2 shows that these are indeed inverse to each other.
\end{exer}

{\bf Definition 3:  The algebraic braid group.} $B_n$ can be regarded
algebraically as the group presented with generators $\sigma_1, \dots,
\sigma_{n-1}$, where $\s_i$ is the braid described in the previous exercise. 
These generators are subject to the relations (as suggested in Figure 2:
$$\sigma_i\sigma_j = \sigma_j\sigma_i, |i-j| > 1,$$
$$\sigma_i\sigma_j\sigma_i = \sigma_j\sigma_i\sigma_j, |i-j| = 1.$$

It was proved by Artin that these are a complete set of relations to abstractly define $B_n$.  This means that any relation among the $\s_j$ can be deduced from the above relations.  

\begin{exer}
Verify that, in $B_n$, all the generators $\s_j$ are conjugate to each other.  Show the total degree
$deg(\b)$ which is the sum of the exponents which appear in an expression of $\b$ in the $\s_j$, is well-defined, i.e. independent of the expression.  Verify that the abelianization of $B_n$, for $n \ge 2$, is 
the infinite cyclic group $\Z$, and $deg: B_n \to \Z$ is equivalent to the abelianization homomorphism.
Thus the commutator subgroup of $B_n$ consists of all words of total degree zero.
\end{exer}

\begin{exer}  
In $B_3$, consider the special element $\Delta = \Delta_3 = \s_1\s_2\s_1$, which consists of a clockwise ``half-twist" of the three strings.  Verify that $\Delta\s_1=\s_2\Delta$ and 
$\Delta\s_2=\s_1\Delta$.  Conclude that $\Delta^2$ commutes with both generators of $B_3$ and therefore is in the centre of $B_3$.  (It actually generates the centre, but we will not prove this.)

Show that in $B_3$ there are two braids $\a$ and $\b$ such that $\a \ne \b$ but $\a^3 = \b^3$.
What about squares?
\end{exer}

\begin{exer}  
For $n \ge 3$, define $\gamma = \gamma_n = (\s_1\s_2\cdots\s_{n-1})^n$.  Draw this braid and convince yourself that it is a ``full-twist'' and that $\gamma$ is central in $B_n$. 
\end{exer}

In fact $\gamma_n$ generates the (infinite cyclic) centre of $B_n$.

We can take a whole countable set of generators $\s_1, \s_2, \dots$ subject to
the above relations, which defines the infinite braid group $B_\infty.$  If we
consider the (non-normal) subgroup generated by $\s_1, \dots, \s_{n-1}$, these
algebraically define $B_n$.  Notice that this convention gives ``natural''
inclusions 
$B_n \subset B_{n+1}$ and $P_n \subset P_{n+1}$.

\begin{exer}  
Does $B_\infty$ have a nontrivial centre?
\end{exer}

Going the other way, if one forgets the last string of an $n+1$-braid the result
is an $n$-braid.  But strictly speaking, this is only a well defined
homomorphism for pure braids, (or at best for the subgroup of braids in which
the string beginning at the point $n+1$ also ends there).  Later, we will have a
use for this forgetful map 
$$f: P_{n+1} \to P_n$$  It is easy to see that $f$ is a left inverse of the
inclusion, or in other words a retraction in the category of groups.

{\bf Artin Combing.}  We now have the ingredients for the combing technique, 
which gives a solution to the word problem for pure braid groups, and
therefore for the full braid groups. 

Solving the word problem in $B_n$ reduces to $P_n$:

Given a word in the $\s_j$, first work
out its corresponding permutation.  If the permutation is nontrivial,
so is the group element it represents.  If trivial, the braid is a
pure braid, and we have reduced the problem to the word problem for $P_n$.

Let $\b$ be a pure $n$-braid and $f(\b)$ the pure $n-1$ braid obtained by
forgetting the last string, but then by inclusion, regard $f(\b)$ in $P_n$. 
Then $\b$ and $f(\b)$ can be visualized as the same braid, except the last string
has been changed in $f(\b)$ so as to have no interaction with the other strings.  Let 
$K$ be the kernel of $f$.  Then it is easy to verify that 
$f(\b)^{-1}\b \in K$ and the map  
$$\b \to (f(\b), f(\b)^{-1}\b)$$ maps $P_n$ bijectively onto the cartesian
product 
$P_{n-1} \times K$.  However, the multiplicative structure is that of a {\it
semidirect} product, as happens whenever we have a split exact
sequence of groups, in this case 
$$1 \to K \to P_n \to P_{n-1} \to 1.$$

 Also notice that every element of $K$ can be represented
by a braid in which the first $n-1$ strands go straight across.  In this way we
identify $K$ with the fundamental group of the complement of the points $\{1,
\dots, n-1\}$ in the plane, which is a free group: $K \cong F_{n-1}$.

\begin{exer}
Draw the braid $\a_i = \s_{n-1}\s_{n-2}\cdots\s_{i+1}\s_i^2\s_{i+1}^{-1}\cdots\s_{n-2}^{-1}\s_{n-1}^{-1}$
and verify that $\a_1, \dots, \a_{n-1}$ represent a free basis for the group $K$ above.
\end{exer}

To solve the word problem in $P_n$: let $\b \in P_n$ be a braid expressed
in the standard generators $\s_i$; the goal is to decide whether $\b$ is really the identity.
Consider its image $f(\b)$, which lies in $P_{n-1}$.  We may assume inductively that the
word problem is solved in $P_{n-1}$.  If $f(\b) \ne 1$ we are done, knowing that $\b \ne 1$.

On the other hand, if $f(\b) = 1$, that means that $\b$ lies in the kernel $K$, which is a free
group.   Express $\b$ as a product of the generators $\a_i$ in $K$, and then we can decide
if $\b$ is the identity, as $K$ is freely generated by the $\a_i$.

The semidirect product decomposition
process can be iterated on $P_{n-1}$ to obtain the Artin normal form: 
$P_n$ is an iterated semidirect product of free groups 
$F_1, F_2, \dots, F_{n-1}$.  Thus every $\b \in P_n$ has a unique expression

$$\b = \b_1\b_2 \cdots \b_{n-1}$$ with $\b_i \in F_i.$

For later reference, we will call the vector
$$(\b_1, \b_2, \cdots, \b_{n-1})$$ the ``Artin coordinates'' of $\b$.

In the above discussion, the subgroup $F_j$ of $P_{j+1}$ (with $j<n$) consists of the pure braids in the subgroup $P_{j+1}$ of $P_n$, which become trivial if the strand labelled $j+1$ is removed.  A set of generators for $F_j$ is
given by the braids $\a_{ij}$, $i<j$, which link the $i$ and $j$ strands over the others:

$$\a_{ij} =  \s_{j-1}\s_{j-2}\cdots\s_{i+1}\s_i^2\s_{i+1}^{-1}\cdots\s_{j-2}^{-1}\s_{j-1}^{-1}$$

Note that $\a_{in}$ is just the braid we called $\a_i$ above.

The semidirect product decomposition permits us to write a finite presentation for $P_n$.
An example is the following:

Generators are $\a_{ij}$ with $1 \le i < j \le n$

Relations are:

$\a_{ij}\a_{ik}\a_{jk} = \a_{ik}\a_{jk}\a_{ij} = \a_{jk}\a_{ij}\a_{ik}$ 
whenever $1 \le i < j < k \le n.$

$\a_{ij}\a_{kl} = \a_{kl}\a_{ij}$ and $\a_{il}\a_{jk} = \a_{jk}\a_{il}$ 
whenever $1 \le i < j < k < l \le n.$

$\a_{ik}\a_{jk}\a_{jl}\a_{jk}^{-1} =  \a_{jk}\a_{jl}\a_{jk}^{-1}\a_{ik}$
whenever $1 \le i < j < k < l \le n.$

\begin{exer}
Verify these relations for specific values of $i, j , k$ and $l$ and convince yourself that they
are ``obvious.''
\end{exer}

\begin{exer}
Show that the abelianization of $P_n$ is free abelian of rank 
$\left(\begin{array}{c}n\\2\end{array}\right)$.  Conclude that the set of generators for $P_n$ given
above is minimal: no smaller set of generators is possible.
\end{exer}

\begin{exer}
By contrast, show that the full braid group $B_n$ can be generated by just two generators.  Hint:
recall that the generators $\s_i$ are all conjugates in $B_n$.
\end{exer}

It may seem strange that even though $P_n$ is a subgroup of $B_n$, its abelianization is, in general, considerably larger.  An even more striking example of this phenomenon can be obtained from the 
free group on two generators, together with its commutator subgroup, which is an infinitely generated free group.  After abelianization they become, respectively, $\Z^2$ and $\Z^\infty$.

We mention an important theorem of W. Thurston.  An infinite group is
called {\it automatic} if is well-modelled (in a well-defined technical sense
which I won't elaborate here) by a finite-state automaton.  The standard
reference is \cite{Ep}.

\begin{thm}
{\bf (Thurston)}  $B_n$ is automatic.
\end{thm}

This implies, for example, that the word problem can be solved by an algorithm
which is quadratic in the length of the input.

\begin{exer}
What is the complexity of the Artin combing algorithm?
\end{exer}

\begin{exer}
Show that $K$ is not normal in $B_n$, but is normalized by the subgroup $B_{n-1}$.  Moreover, the action of $B_{n-1}$ on $K$ by conjugation is essentially the same as the Artin presentation, given in the next section. 
\end{exer}

\section{Mapping class groups and braids}

{\bf Definition 4:  $B_n$ as a mapping class group.}  Going back to the second
definition, imagine the particles are in a sort of planar jello and pull their
surroundings with them as they dance about.  Topologically speaking, the motion
of the particles extends to a continuous family of homeomorphisms of the plane
(or of a disk, fixed on the boundary).  This describes an equivalence between
$B_n$ and the mapping class of $D_n$, the disk $D$ with $n$ punctures (marked
points).  That is, $B_n$ can be considered as the group of homeomorphisms of $D_n$ 
fixing $\partial D$ and permuting the punctures,
modulo isotopy fixing $\partial D \cup \{1, \dots, n\}$.

{\bf Definition 5:  $B_n$ as a group of automorphisms.}  A mapping class $[h]$,
where $h:D_n \to D_n$, gives rise to an automorphism $h_*: F_n \to F_n$ of free
groups.  Using the interpretation of braids as mapping classes, this
defines a homomorphism

$$B_n \to Aut(F_n),$$

which Artin showed to be faithful, i. e. injective.

\begin{center}%
\includegraphics[totalheight=5cm]{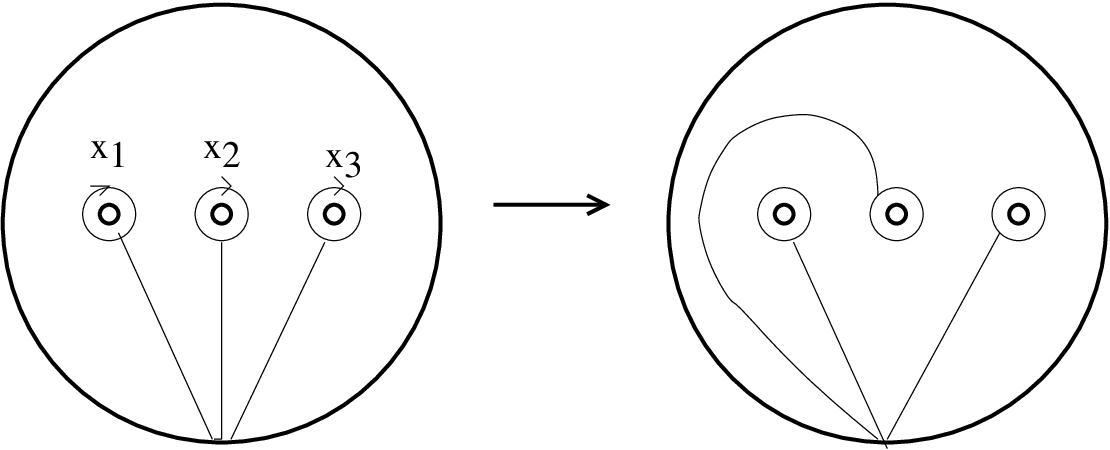}%
\begin{center}%
Figure 3: The action of $\s_1$ on $D_3$.
\end{center}
\end{center}

The generator $\sigma_i$ acts as

\begin{align}
x_i &\to x_i x_{i+1}x_i^{-1} \cr x_{i+1} &\to x_i \cr x_j &\to x_j,
j \ne i, i+1. 
\end{align}

This is shown in the figure above, for the case $n=3$, and the action of $\s_1$.

\begin{thm}
\label{artin}
{\bf (Artin)} Under the identification described above, $B_n$ is
the set of automorphisms $h \in Aut(F_n)$ of the form
$$h(x_i) =  w_j^{-1}x_jw_j,$$ where $w_j$ are words in $F_n$, and satisfying
$h(x_1 \cdots x_n) = x_1 \cdots x_n$.
\end{thm}
  
This point of view gives further insight into the group-theoretic properties of
braid groups.  A group $G$ is {\it residually finite} if for every 
$g\in G$ there is a homomorphism $h:G \to F$ onto a finite group $F$ such that
$h(g)$ is not the identity.   In other words, any element other than the identity can be proved nontrivial by looking at some homomorphism of $G$ to a finite group.

\begin{exer}
Show that subgroups of residually finite groups are residually finite.  
\end{exer}

\begin{exer}
Show that the automorphism group of any finitely generated residually finite group is itself residually finite.  [This is a result of Baumslag.  You may want to consult \cite{KMS} for a proof, and further discussion on the subject.]
\end{exer}

\begin{exer}  Show that finitely-generated abelian groups are residually finite.
The matrices $\left[\begin{array}{cc}1 & 2 \\0 & 1\end{array}\right]$ and
$\left[\begin{array}{cc}1 & 0 \\2 & 1\end{array}\right]$ generate a free group in the 
group $SL(2,\Z)$, two-by-two integral matrices with determinant $1$, which is
a group of automorphisms of $\Z^2$.  Thus $SL(2,\Z)$ is residually finite and so is the rank 2 free group.  Conclude that a free group of any finite or countable rank is residually finite.
\end{exer}

From these observations, we see that $Aut(F_n)$ is residually finite and conclude the following.

\begin{thm}
$B_n$ is residually finite.
\end{thm}

A group $G$ is said to be Hopfian if it is not isomorphic with any nontrivial quotient group.  In other words, if $G \to H$ is any surjective homomorphism with nontrivial kernel, then $H \not\cong G$.

\begin{exer}
A finitely-generated residually finite group is Hopfian.
\end{exer}

It follows immediately that

\begin{thm}
$B_n$ is Hopfian.
\end{thm}

{\bf The word problem:}  If a group $G$ is given by generators, say $\{g_i\}$, and relations, it may be difficult to determine whether two words in the generators actually represent the same element of $G$.  This is called the {\em word problem}.  It easily reduces to the question: given a word in the $g_i$, does it present the identity element of $G$? 

If $G$ is residually finite, then in principle a give word which does {\it not} represent the identity, can be proven so by finding a homomorphism of $G$ onto a finite group in which the word is sent to a nontrivial element, where verifying this is a finite problem.  Of course, finding such a homomorphism may well be terribly difficult.  But for particular groups there may be more algorithmic methods for solving the word problem.  Giving such an algorithm, for a particular group, or class of groups, is what is meant by a 
{\em solution} to the word problem.

Such solutions exist for free groups and for the braid groups.  First, consider a free group, with free generators $\{x_i\}$.  If $w$ is a given word in these generators, simply perform all free cancellations which are possible: that is, remove any two consecutive letters which happen to be $x_ix_i^{-1}$ or 
$x_1^{-1}x_i$.  Continue doing this until no such cancellations are possible.  If the result is the empty word, $w$ represents the identity; otherwise it does not.

Artin's theorem provides a means for solving the word problem algorithmically in $B_n$ by reducing it to several word problems in a free group.  Suppose a word $w$ in the Artin generators $\s_i$ is given.
Then we may consider $w$ to be an automorphism of the free group $F_n$.   Using equations (1), we can explicitly calculate the values $w(x_1), \dots, w(x_n)$ as words in the $x_i$.  We can then solve the word problem $n$ times in $F_n$ to decide whether $w(x_i) = x_i$ for each 
$i = 1, \dots , n$.  If this is the case, $w$ represents the identity of $B_n$; otherwise, it does not.  

Another solution to the word problem in braid groups, again reducing it to free groups, will be discussed later, using the technique called ``Artin combing.''  Yet another, due to Garside, also solves the so-called
{\em conjugacy} problem, which is to decide whether two given words in the generators are conjugates in the group.  I will not discuss Garside's method here, as it is covered quite thoroughly in the notes on cryptography in this workshop.

{\bf The modular group.}  There is an interesting connection between the braid group $B_3$ and 
two-by-two integral matrices. 

Consider the matrices $$S = \left[\begin{array}{cc}0 & -1 \\1 & 0\end{array}\right] \quad\quad 
T = \left[\begin{array}{cc} 0 & -1 \\ 1 & 1\end{array}\right] .$$
One can easily check that $S^2 = -I = T^3 $.  These matrices generate cyclic subgroups of $SL(2,\Z)$,
$\Z_4$ and $\Z_6$, respectively (we use the abbrviation $\Z_n = \Z/(n\Z)$).  It is well-known that, using the generators $S$ and $T$,  there is an amalgamated free product structure
$$SL(2,\Z) \cong \Z_4 *_{\Z_2} \Z_6.$$

The {\em modular group} $PSL(2,\Z)$ is the quotient of $SL(2,\Z)$ by its center $\{\pm I\}$.  
One can also regard $PSL(2,\Z)$ as the group of M\"obius transformations of the complex plane 
$$z \mapsto {\frac{az + b}{cz + d}} \quad{\rm corresponding \: to} \quad
\left[\begin{array}{cc}a & b \\c & d\end{array}\right] \in SL(2,\Z).$$

The modular group has the structure of a free product $\Z_2*\Z_3$, with (the cosets of) $S$ and $T$ generating the respective factors.

\begin{exer}
Show that the braid group $B_3 \cong  \langle \s_1, \s_2 | \s_1\s_2\s_1 =  \s_2\s_1\s_2 \rangle$ also has the presentation $\langle x, y | x^2 = y^3 \rangle$, by finding appropriate expressions of $x, y$ in terms of $\s_1, \s_2$, and vice-versa, so that the transformations respect the relations and are mutual inverses.

Verify that the element $x^2 (=y^3)$ is central in $B_3$ and in fact generates the center of $B_3$.

Show that $B_3$ modulo its center is isomorphic with $PSL(2,\Z)$.

\end{exer}

{\bf Definition 6: $B_n$ as a fundamental group.}  In complex $n$-space $\C^n$
consider the big diagonal 
$$\Delta = \{(z_1, \dots, z_n) ;\quad z_i = z_j, {\quad \rm some \quad} i < j \}
\subset
\C^n.$$ 

Using the basepoint $(1,2, \dots, n)$, we see that 
$$P_n = \pi_1(\C^n \setminus \Delta).$$

In other words, pure braid groups are fundamental groups of complements of a
special sort of complex {\it hyperplane arrangement}, itself a deep and
complicated subject.

To get the full braid group we need to take the fundamental group of the {\it
configuration space}, of orbits of the obvious action of $\Sigma_n$ upon $\C^n
\setminus \Delta$. Thus
$$B_n = \pi_1((\C^n \setminus \Delta) / \Sigma_n).$$

Notice that since the singularities have been removed, the projection
$$\C^n \setminus \Delta  \longrightarrow
(\C^n \setminus \Delta) / \Sigma_n$$ is actually a covering map.  As is
well-known, covering maps induce injective homomorphisms at the $\pi_1$ level,
so this is another way to think of the inclusion $P_n \subset B_n$.

It was observed in \cite{FN} that $\C^n \setminus \Delta$ has trivial homotopy
groups in dimension greater than one.  That is, it is an Eilenberg-Maclane space,
also known as
a  $K(P_n, 1)$.  Therefore its cohomology groups coincide with the group cohomology
of $P_n$.  By covering theory, the quotient space
$(\C^n \setminus \Delta) / \Sigma_n$ also has trivial higher homotopy, so it is a
$K(B_n, 1).$  Since these spaces have real dimension $2n$, this view of braid
groups gives us the following observation.

\begin{thm} 
The groups $B_n$ and $P_n$ have finite cohomological dimension.
\end{thm}

If a group contains an element of finite order, standard homological algebra
implies that the cohomological dimension of the group must be infinite.  Thus
there are no braids of finite order.  

\begin{cor} 
The braid groups are torsion-free.
\end{cor}

Finally, we note that the space 
$(\C^n \setminus \Delta) / \Sigma_n$ can be identified with the space of all
complex polynomials of degree $n$ which are monic and have $n$ distinct roots

$$p(z) = (z - r_1) \cdots (z - r_n).$$

This is one way in which the braid groups play a role in classical algebraic
geometry, as fundamental groups of such spaces of polynomials.

{\bf Definition 7:}  $B_n$ is the fundamental group of the space of all monic polynomials
of degree $n$ with complex coefficients and only simple roots.

\subsection*{Surface braids}

Let $\Sigma$ denote a surface, with or without boundary.  One can define a braid on $\Sigma$
in exactly the same way as on a disk -- namely, a collection of disjoint paths $\b_1(t), \cdots, \b_n(t)$
in $\Sigma \times I$, with each $\b_i(t) \in \Sigma \times \{t\}$ so that they begin and end at the same points of $\Sigma$, possibly permuted.  They form a group, which is usually denoted $B_n(\Sigma)$.

The subgroup of braids for which the permutation of the points is the identity is the {\it pure} surface braid group $P_n(\Sigma)$.  Thus $B_n(D^2)$ and $P_n(D^2)$ are the classical Artin braid groups which we
have discussed before.  Much of the above discussion holds for surface braid groups; for example, they may be interpreted as fundamental groups of configuration spaces.

On the other hand, there are differences.  For example the braid groups of the sphere $S^2$ and the
projective plane $\R P^2$ have elements of finite order.  Van Buskirk showed that these are the only closed surfaces to have torsion in their braid groups, however.  There are finite presentations for surface braid groups in the literature, and these groups are currently an active area of research.  One reason for their interest is that they play an important part in certain topological quantum field theories (TQFT's).

\section{Knot theory, braids and the Jones polynomial}

\subsection*{Knots and Reidemeister moves}

One often pictures a knot or link by drawing a projection onto the plane, with only
double points, and indicating which string goes under by putting a small gap in it, as in
Figure 1.  If one deforms the plane by a planar isotopy, an equivalent projection results.
One can also make local changes to a knot diagram.  The first two vignettes in Figure 2 are
sometimes called Reidemeister moves of type 2 and 3 (respectively).  One can remember
the numbering, as type 2 involves two strands and type 3 involves three.  There is also a
type 1 Reidemeister move.

\begin{center}%
\includegraphics[totalheight=2cm]{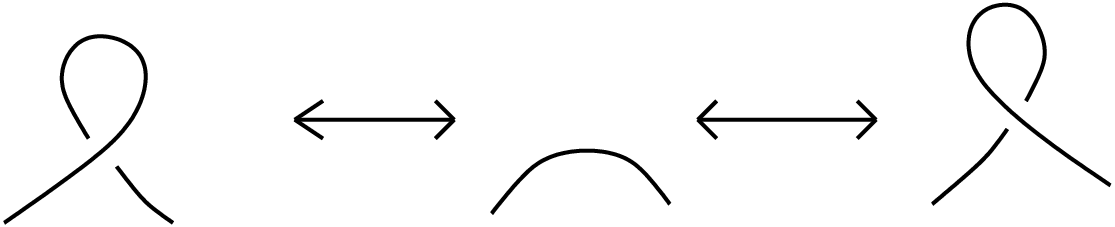}%
\begin{center}%
Figure 4: Reidemeister move of type 1
\end{center}
\end{center}

\begin{thm}[Reidemeister] 
Two knot diagrams present equivalent knots if and only if one can
be transformed to the other by a finite sequence of Reidemeister moves (plus an isotopy of the plane).
\end{thm}

\subsection*{Markov moves}
It has already been mentioned that different braids can, upon forming
the closure, give rise to the same knot or link.  You can easily
convince yourself that if $\a$ and $\b$ are $n$-string braids, then
the closure of $\a^{-1}\b\a$ is equivalent to the closure of $\b$.
Because of the trivial strings added in forming the closure, $\a$ and
its inverse can annihilate each other!
Therefore we see that conjugate braids close to the same link.  

Another example of
this phenomenon is to consider the $n$-braid $\b$ as an element of
$B_{n+1}$ by adding a trivial string at the top,
 and compare the original closure $\hat{\b}$ with the closure 
of $\b \times \s_n$, the latter braid taken in $B_{n+1}$. If you
sketch this, you can easily convince yourself that these result in the
same link.  The first move (conjugation) and the second move just
described constitute the two Markov moves.  They are the key to
understanding the connection between braid theory and knot theory.
The following theorem is due to J. W. Alexander (first part) and
A. A. Markov.  A full proof appeared first, to my knowledge, in
\cite{Bir1}.

\begin{thm}
Every knot or link is the closure of some braid.  Two braids close to
equivalent knots or links if and only if they are related by a finite
sequence of the two Markov moves.
\end{thm}

\begin{exer}
If you want to learn a nice, elementary proof of the above theorem, read 
the paper \cite{morton}.
\end{exer}

\begin{exer}
Consider the 3-braid $\b = \s_1\s_2$ and identify the knot or link, simplifying, if possible
by using Reidemeister moves: $\hat{\b}$, $\hat{\b^2}$, $\hat{\b^3}$.  For which integers $k$
is the closure of $\b^k$ a knot?
\end{exer}

\subsection*{Kauffman's bracket and Jones' polynomial}
The original construction of the Jones polynomial involved a family of representations
$$B_n \to {\mathcal A}$$ of the braid groups into an algebra ${\mathcal A}$ involving a parameter $t$.
This algebra has a linear trace function into the ring of (Laurent) polynomials.  The composite of the representation and the trace, together with some correction terms to account for the Markov moves, defined the original Jones polynomial.

We will present this in a sort of reverse order, as there is a very simple derivation of the Jones polynomial discovered a few years later by L. Kauffman.  From this, we can define an algebra, called the Temperley-Lieb algebra, and reconstruct what is essentially Jones' representation.

Consider a planar diagram $D$ of a link $L$, which has only simple transverse crossings.  Kauffman's bracket $\langle D \rangle$ is (at first) a polynomial in three commuting variables, $a, b$ and $d$ defined by the equations:

\begin{center}%
\includegraphics[totalheight=3cm]{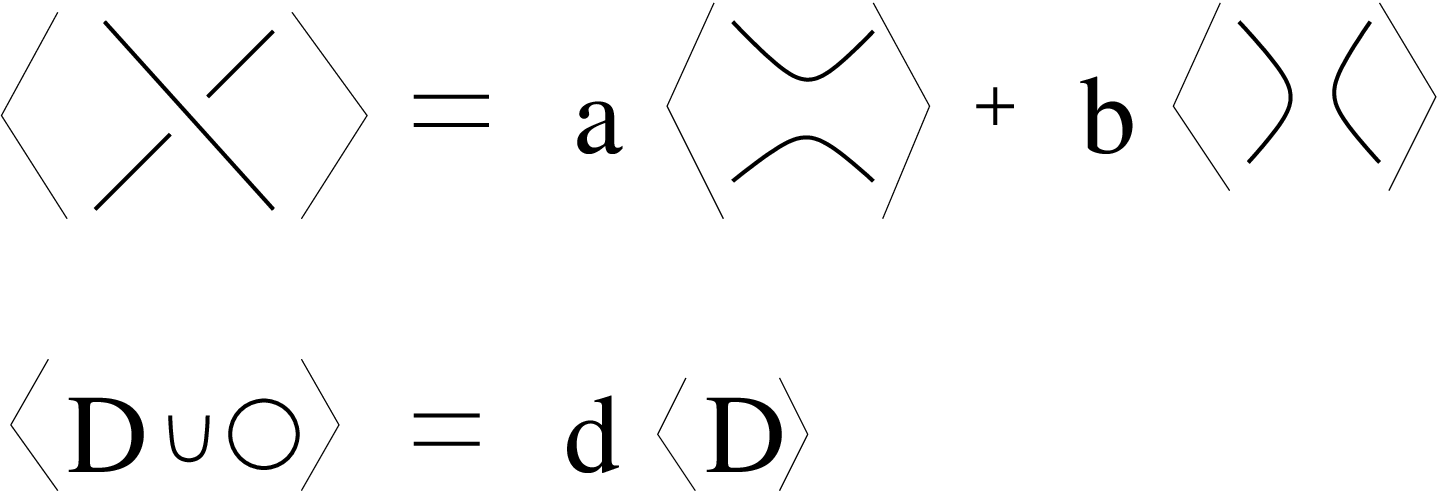}%
\begin{center}%
Figure 5: Equations defining Kauffman's bracket polynomial.
\end{center}
\end{center}

Some explanation is in order.  The vignettes in the brackets in the first equation stand for complete link diagrams, which differ only near the crossing in question.  Those on the right side represent diagrams with fewer crossings.  In the second equation, one can introduce a closed curve which has no intersections with the remainder of the diagram $D$, resulting in a diagram whose bracket polynomial is $d$ times the bracket of $D$.  

\begin{exer}
Verify that the bracket polynomial is well-defined, if we decree that the bracket of a single curve with no crossings is equal to 1.
\end{exer}

\begin{exer}
Show that  the bracket will be invariant under the type 2 Reidemeister move, if we have the following relations among the variables:
$$a^2 + b^2 + abd = 0 \quad {\rm and } \quad ab = 1.$$
\end{exer}

Thus we make the substitutions $b = a^{-1}$ and $d = - a^2 - a^{-2}$ and now consider the bracket to be a Laurent polynomial in the single variable $a$. 

\begin{exer}
Show that  the  invariance under the type 2 Reidemeister moves implies the bracket is invariant under the type 3 move, too.  
\end{exer}

\begin{exer}
Calculate the bracket of the two trefoil knot diagrams with 3 crossings.  Show that they are the same, except for reversal of the sign of the exponents.
\end{exer}

\begin{exer}
Investigate the effect of Reidemeister move 1 on the bracket of a diagram, and show that it changes the bracket by a factor of $-a^{\pm 3}$, the sign of the exponent depending on the sense of the curl removed.
\end{exer}

Because of this, one can define a polynomial invariant under all three Reidemeister moves, by counting the number of positive minus the number of negative crossings, and modifying the bracket polynomial by an appropriate factor.  A positive crossing corresponds to a (positive) braid generator, if both strings are oriented from left to right.  This gives, up to change of variable, the Jones polynomial of the knot.
Specifically, we define the writhe of an {\it oriented} diagram $\vec D$ for an oriented knot (or link) $\vec K$ to be
$$w(D) = \sum_c  \epsilon_c$$
where the sum is over all crossings and $\epsilon_c = 1$ if the crossing $c$ is positive, and $-1$ if negative.  Then we define
$$f_{\vec K}(a) = (-a^3)^{-w(\vec D)} \langle D \rangle.$$
   This is an invariant of the oriented link $\vec K$.  If $K$ happens to be a knot, it is independent of the orientation, as reorienting both strands of a crossing does not change its sign.  It is related to the Jones polynomial $V_K(t)$ by a simple change of variables:
 $$V_{\vec K}(t) = f_{\vec K}(t^{-1/4}).$$ 
 
 \begin{exer}
Show that all exponents of $a$ in $\langle D \rangle$ are divisible by 4 if $D$ is a diagram of a knot (or link with an odd number of components).  In other cases, they are congruent to 2 mod 4.  Thus the Jones polynomial is truly a (Laurent) polynomial in $t$ for knots and odd component links, but a polynomial in $\sqrt{t}$ in the other cases.
\end{exer}

\section*{Representations}

This is a very big subject, which I will just touch upon.  By a representation of
a group we will mean a homomorphism of the group into a group of matrices, or
more generally into some other group, or ring or algebra.  Often, but not always,
we want the target to be finite-dimensional.  We've already encountered the Artin
representation $B_n \to Aut(F_n)$, which is faithful.  Here the target group is
far from being ``finite-dimensional.''

Another very important representation is the one defined by Jones \cite{J} which gave
rise to his famous knot polynomial, and the subsequent revolution in knot theory. 
The version I will discuss is more thoroughly described in \cite{K}; it is based on
the Kauffman bracket, an elementary combinatorial approach to the Jones polynomial.
First we need to describe the Temperley-Lieb Algebras $\T_n$, in their geometric
form.  The elements of $\T_n$ are something like braids: we consider strings in a
box, visualized as a square in the plane, endpoints being exactly $n$ specified
points on each of the left and right sides.   The strings are not required to be
monotone, or even to run across from one side to the other.  There also may be
closed components.  Really what we are looking at are ``tangle'' diagrams.  Two
tangle diagrams are considered equal if there is a planar isotopy, fixed on the
boundary of the square, taking one to the other.

\begin{center}%
\includegraphics[totalheight=3cm]{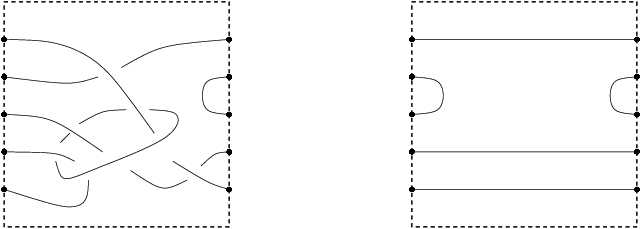}%
\begin{center}%
Figure 6: A typical element in $\T_5$ and the generator $e_3$.
\end{center}
\end{center}


Now we let $A$ be a fixed complex
number (regarded as a parameter), and formally define $\T_n$ to be the complex
vector space with basis the set of all tangles, as described above, but modulo
the following relations, which correspond to similar
relations used to define Kauffman's bracket version of the Jones polynomial (we have promoted the variable $a$ to upper case).

\begin{center}%
\includegraphics[totalheight=4cm]{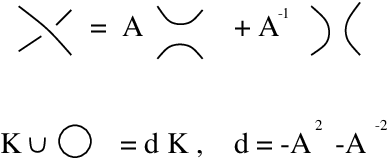}%
\begin{center}%
Figure 7: Relations in $\T_n$.
\end{center}
\end{center}


The first relation means that we can replace a tangle with a crossing by a linear combination of two tangles with that crossing removed in two ways.  As usual, the pictures mean that the tangles are identical outside the part pictured.  The second relation means that we can remove any closed curve in the diagram, if it does not have any crossings with the rest of the tangle, at the cost of multiplying the tangle by the scalar $-A^2 -A^{-2}$.
Using the relations, we see that any element of $\T_n$ can be expressed as a
linear combination of tangles which have no crossings and no closed curves --
that is, disjoint planar arcs connecting the $2n$ points of the boundary.  This
gives a finite generating set, which (for generic values of $A$) can be shown to
be a basis for $\T_n$ as a vector space.  But there is also a multiplication of
$\T_n$, a concatenation of tangles, in exactly the same way braids are multiplied.
This enables us to consider $\T_n$ to be generated {\it as an algebra} by the
elements $e_1, \dots, e_{n-1}$.  In $e_i$ all the strings go straight across, except those at level $i$ and $i+1$ which are connected by short caps; the generator $e_3$ of $\T_5$ is illustrated in Figure 2.  The identity of this algebra is simply the
diagram consisting of $n$ horizontal lines (just like the identity braid).  $\T_n$ can be described abstractly as the associative algebra with the generators
$e_1, \dots, e_{n-1}$, subject to the relations:

$$
e_ie_j = e_je_i \;{\rm if}\; |i-j|>1\quad\quad
e_ie_{i\pm 1}e_i = e_i \quad\quad
e_i^2 = (-A^2 - A^{-2})e_i
$$

It is an enjoyable exercise to verify these relations from the pictures.  Now the
Jones representation $J: B_n \to \T_n$ can be described simply by considering a
braid diagram as an element of the algebra.  In terms of generators, this is just

$$
J(\s_i) = A + A^{-1}e_i
$$

\subsection*{The Burau representation}

One of the classical representations of the braid groups is the Burau
representation, which can be described as follows.  Consider the definition of
$B_n$ as the mapping class group of the punctured disk $D_n$ (Definition 4).
As already noted, the fundamental group of $D_n$ is a free group, with generator
$x_i$ represented by a loop, based at a point on the boundary of the disk, which
goes once around the $i^{th}$ puncture.  Consider the subgroup of $\pi_1(D_n)$
consisting of all words in the $x_i$ whose exponent sum is zero.  This is a normal
subgroup, and so defines a regular covering space 
$$p \colon \bar{D}_n \to D_n.$$ The group
of covering translations is infinite cyclic.  Therefore, the homology
$H_1(\bar{D}_n)$ can be considered as a module over the polynomial ring
$\Lambda := \Z[t,t^{-1}]$, where $t$ represents the generator of the covering translation
group.  Similarly, if $*$ is a basepoint of $D_n$, the relative homology group $H_1(\bar{D}_n, p^{-1}(*))$ is a $\Lambda$-module.

 A braid $\b$ can be represented as (an isotopy class of) a homeomorphism
$\b :D_n \to D_n$ fixing the basepoint.  This lifts to a homeomorphism
$\bar{\b} : \bar{D_n} \to \bar{D_n}$, which is unique if we insist that it fix
some particular lift of the basepoint.  The induced homomorphism on homology,
$\bar{\b}_* : H_1(\bar{D_n},p^{-1}(*)) \to H_1(\bar{D_n},p^{-1}(*))$ is a linear map of these
finite-dimensional modules, and so can be represented by a matrix with entries
in $\Z [t,t^{-1}]$.  The mapping
$$\b \to \bar{\b}_*$$ is the Burau representation of $B_n$.

Let's illustrate this for the case $n=3$.  $D_3$ is replaced by the wedge of three
circles, which is homotopy equivalent to it, to simplify visualization.  The covering space $\bar{D_n}$
is shown as an infinite graph.  Although as an abelian group $H_1(\bar{D_3},p^{-1}(*))$ is infinitely generated,
as a $\Lambda$-module it has three generators $g_i = \tilde{x_i}$, the lifts of the generators $x_i$ of
$D_3$, $i = 1, 2, 3$ beginning at some fixed basepoint in $p^{-1}(*)$.  The other elements of  
$H_1(\bar{D_3},p^{-1}(*))$ are represented by translates $t^ng_i$.

Now consider the action of $\s_1$ on $D_3$, at the fundamental group level.  As we saw in Definition
5, $\s_{1*}$ takes $x_1$ to $x_1x_2x_1^{-1}$.  Accordingly $\bar{\s_1}_*$ takes $g_1$ to the lift of $x_1x_2x_1^{-1}$,
which is pictured in the lower part of the illustration.  In terms of homology, this is $g_1 + tg_2 -tg_1$.  
Therefore $$\bar{\s_1}_*(g_1) = (1-t)g_1 + t g_2.$$
Similarly$$\bar{\s_1}_*(g_2) = g_1, \quad \bar{\s_1}_*(g_3) = g_3.$$

\begin{center}%
\includegraphics[totalheight=5cm]{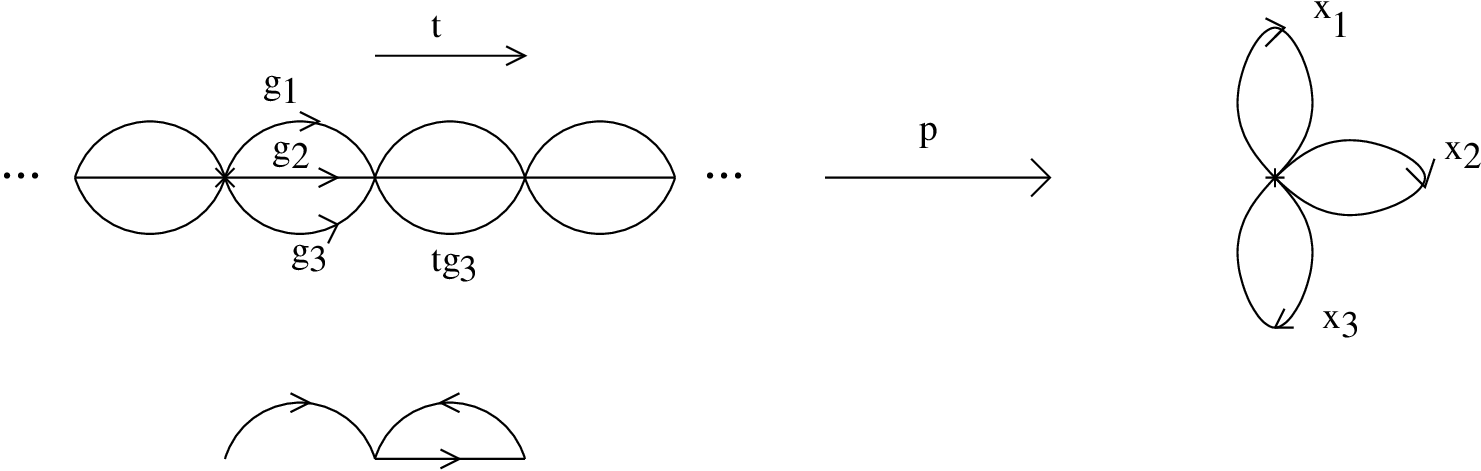}%
\begin{center}%
Figure 8:  Illustrating the Burau representation for $n=3$.
\end{center}
\end{center}

\begin{exer}
Show that, with appropriate choice of basis, the Burau representation of $B_n$ 
sends $\s_j$ to the matrix
 $$I_{j-1} \oplus \left[\begin{array}{cc}1-t & t \\1 & 0 \end{array}\right] \oplus I_{n-j-1},$$ 
 where $I_k$ denotes the $k \times k$ identity matrix.
\end{exer}

{\bf A probabilistic interpretation:}  Vaughan Jones offered the following interpretation of the Burau representation, which I will modify slightly because we have chosen opposite crossing conventions.  Picture a braid as a system of trails.  Wherever one trail crosses under another there is a probability $t$ that a person will jump from the lower trail to the upper trail; the probability of staying on the same trail 
is $1 - t$.  The $i,j$ entry of the Burau matrix corresponding to a braid then represents the probability that, if a person starts on the trail at level $i$, she will finish on the trail at level $j$.

It has been known for many years that this representation is faithful for 
$n \le 3$, and it is only within the last decade that it was found to be
unfaithful for any $n$ at all.  John Moody \cite{Mo} showed in 1993 that it is unfaithful
for $n \ge 9$.  This has since been improved by Long and Paton \cite{LP} and more recently by Bigelow to $n \ge 5$. The case $n = 4$ remains open, as far as I am aware.

\subsection*{Linearity of the braid groups.}  

It has long been questioned, whether
the braid groups are {\it linear}, meaning that there is a faithful
representation  $B_n \to GL(V)$ for some finite dimensional vector space $V$.  A
candidate had been the so-called Burau representation, but as already mentioned it
has been known for several years that Burau is unfaithful, in general. 
  
The question was finally settled recently by Daan Krammer and Stephen
Bigelow, using equivalent representations, but different methods.
They use a representation defined very much like the Burau representation.  But
instead of a covering of the punctured disk $D_n$, they use a covering of the
{\it configuration space} of pairs of points of $D_n$, upon which $B_n$ also acts.
This action induces a linear representation in the homology of an appropriate
covering, and provides just enough extra information to give a faithful
representation!

\begin{thm}
{\bf (Krammer \cite{krammer}, Bigelow \cite{bigelow})} The braid groups are linear. 
\end{thm}

In fact, Bigelow has announced that the BMW representation (Birman, Murakami,
Wenzl) \cite{BW} is also faithful.   Another open question is whether the Jones
representation $J: B_n \to \T_n$, discussed earlier,
is faithful. 

\section{Ordering braid groups}  

This is, to me, one of the most exciting of the recent developments in braid theory. 
Call a group $G$ {\it right orderable} if its elements can be given a strict
total ordering $<$ which is right-invariant: 
$$
\forall x,y,z \in G, \quad x<y \Rightarrow xz < yz.
$$

\begin{thm}
{\bf (Dehornoy\cite{De})}  $B_n$ is right-orderable.
\end{thm}

Interestingly, I know of three quite different proofs.  The first
is Dehornoy's, the second is one that was discovered jointly by myself and four
other topologists.  We were trying to understand difficult technical aspects of
Dehornoy's argument, then came up with quite a different way of looking at exactly
the same ordering, but using the view of $B_n$ as the mapping class group of the
punctured disk $D_n$.  Yet a third way is due to Thurston, using the fact that
the universal cover of $D_n$ embeds in the hyperbolic plane.  Here are further
details.

Dehornoy's approach: It is routine to verify that a group $G$ is right-orderable
if and only if there exists a subset
$\Pi$ (positive cone) of $G$ satisfying:

(1) $\Pi\cdot\Pi = \Pi$

(2)  The identity element does not belong to $\Pi$, and for every $g \ne 1$ in $G$
exactly one of $g \in \Pi$ or $g^{-1} \in \Pi$ holds.

One defines the ordering by $g<h$ iff $hg^{-1} \in \Pi$.  

\begin{exer}
Verify that the transitivity law holds, and that the ordering is right-invariant.  Show that the ordering
defined by this rule is also left-invariant if and only if $g\Pi g^{-1} = \Pi$ for every $g \in G$; that is,
$\Pi$ is ``normal.''
\end{exer}

Dehornoy's idea is to call a braid $i$-positive if it is expressible as a word
in $\s_j, j\ge i$ in such a manner that all the exponents of $\s_i$ are
positive.  Then define the set $\Pi \subset B_n$ to be all braids which are
$i$-positive for some $i = 1, \dots, n-1.$  To prove (1) above is quite easy,
but (2) requires an extremely tricky argument.

Here is the point of view advocated in \cite{FGRRW}.  Consider $B_n$ as acting on the
complex plane, as described above. Our idea is to consider the image of the real
axis $\b(\R)$, under a mapping class $\b \in B_n$.  Of course there are choices
here, but there is a unique ``canonical form'' in which (roughly speaking) $\R
\cap \b(\R)$ has the fewest number of components.  Now declare a braid $\b$ to be positive
if (going from left to right) the first departure of the canonical curve
$\b(\R)$ from $\R$ itself is into the {\it upper} half of the complex plane. 
Amazingly, this simple idea works, and gives exactly the same ordering as 
Dehornoy's combinatorial definition.

Finally, Thurston's idea for ordering $B_n$ again uses the mapping class point
of view, but a different way at looking at ordering a group.  This approach, which has the advantage of defining infinitely many right-orderings of $B_n$ is described by H. Short and B. Wiest in \cite{SW}.  The Dehornoy ordering (which is discrete) occurs as one of these
right-orderings -- others constructed in this way are order-dense.  A group $G$ acts
on a set $X$ (on the right) if the mapping $x \to xg$ satisfies:
$x(gh) = (xg)h$ and $x1 = x$.  An action is {\it effective} if the only element
of $G$ which acts as the identity is the identity $1 \in G$.  The following is a
useful criterion for right-orderability:

\begin{lemma}
If the group $G$ acts effectively on $\R$ by order-preserving
homeomorphisms, then $G$ is right-orderable.
\end{lemma}

By way of a proof, consider a well-ordering of the real numbers.  Define, for
$g$ and $h \in G$,

$$g < h \Leftrightarrow xg < xh {\rm\quad at\ the\ first\ } x \in \R {\rm\ such\
that\ } xg \ne xh.$$

It is routine to verify that this defines a right-invariant strict total
ordering of $G$.  (By the way, we could have used any ordered set in place of
$\R$.)  For those wishing to avoid the axiom of choice (well-ordering $\R$) we could have just used an ordering of the rational numbers.

The universal cover $\tilde{D_n}$ of $D_n$ can be
embedded in the hyperbolic plane $\H^2$ in such a way that the covering
translations are isometries. This gives a hyperbolic structure on $D_n$.  It also
gives a beautiful tiling of $\H^2$, illustrated in Figure 9 for the case $n=2.$

\begin{center}%
\includegraphics[totalheight=10cm]{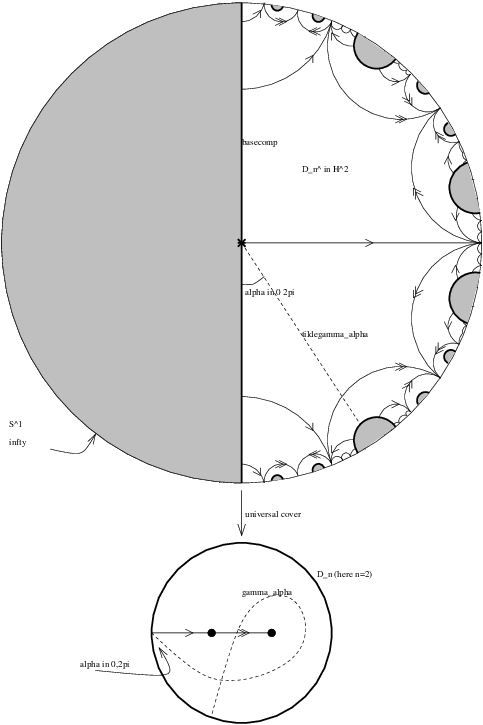}%
\begin{center}%
Figure 9: The universal cover of a twice-punctured disk, with a lifted geodesic. 
(Courtesy of H. Short and B. Wiest \cite{SW}
\end{center}
\end{center}


Choose a
basepoint $* \in \partial D$ and a specific lift
$\tilde{*} \in \H^2$.  Now a braid is represented by a homeomorphism of $D_n$,
which fixes $*$.  This homeomorphism lifts to a homeomorphism of $\tilde{D_n}$,
unique if we specify that it fixes $\tilde{*}$.  In turn, this homeomorphism
extends to a homeomorphism of the boundary of $\tilde{D_n}$  But in fact,
this homeomorphism fixes the interval of $\partial\tilde{D_n}$ containing 
$\tilde{*}$, and if
we identify the complement of this interval with the real line $\R$, a braid
defines a homeomorphism of $\R$.  This defines an action of $B_n$ upon $\R$ by
order-preserving homeomorphisms, and hence a right-invariant ordering of $B_n$.

{\bf Two-sided invariance?}  Any right-invariant ordering of a group can be
converted to a left-invariant ordering, by comparing inverses of elements, but
that ordering is in general different from the given one.  We will say that a group
$G$ with strict total ordering $<$ is fully-ordered, or bi-ordered, if 
$$x<y \Rightarrow xz<yz {\rm\ and\ } zx<zy, \forall x,y,z \in G.$$

There are groups which are right-orderable but not bi-orderable -- in fact
the braid groups!

\begin{prop}
{\bf (N. Smythe)} For $n>2$ the braid group $B_n$ cannot be
bi-ordered.
\end{prop}

The reason for this is that there exists a nontrivial element which is conjugate
to its inverse: take
$x=\s_1\s_2^{-1}$ and $y = \s_1\s_2\s_1$ and note that 
$yxy^{-1} = x^{-1}.$  In a bi-ordered group, if $1<x$ then $1<yxy^{-1}=x^{-1}$,
contradicting the other conclusion
$x^{-1} < 1$.  If $x<1$ a similar contradiction arises.

\begin{exer}
Show that in a bi-ordered group $g < h$ and $g' < h'$ imply $gg' < hh'$.  Conclude that if 
$g^n = h^n$ for some $n \ne 0$, then $g = h$.  That is, roots are unique.  Use this to give an alternative proof that $B_n$ is not bi-orderable if $n \ge 3$.
\end{exer}

\begin{thm}
The pure braid groups $P_n$ can be bi-ordered.
\end{thm}

This theorem was first noticed by J. Zhu, and the argument
appears in \cite{RZ}, based on the result of Falk and Randall \cite{FR} that the pure
braid groups are ``residually torsion-free nilpotent.''  Later, in joint work
with  Djun Kim, we discovered a really natural, and I
think beautiful, way to define a bi-invariant ordering of $P_n$.  We've already
done half the work, by discussing Artin combing.  Now we need to discuss
ordering of free groups.

\subsection*{Bi-ordering free groups}

\begin{lemma}\label{orderfree}
For each $n \ge 1$, the free group $F_n$ has a bi-invariant
ordering $<$ with the further property that it is invariant under any
automorphism
$\phi :F_n \to F_n$ which induces the identity upon abelianization: $\phi_{ab} =
id :\Z^n \to \Z^n$. 
\end{lemma}

The construction depends on the Magnus expansion of free groups into rings of
formal power series.  Let 
$F$ be a free group with free generators $x_1, \dots, x_n$. Let
$Z[[X_1, \dots,X_n]]$ denote the ring of formal power series in the {\it
non-commuting} variables
$X_1, \dots,X_n$.  

Each term in a formal power series has a well-defined (total)
degree, and we use $O(d)$ to denote terms of degree $\ge d.$  The subset $\{1 +
O(1)\}$ is actually a multiplicative subgroup of 
$Z[[X_1, \dots,X_n]]$.  Moreover, there is a multiplicative homomorphism
$$\mu: F \to Z[[X_1, \dots,X_n]]$$ defined by

\begin{align}
\mu(x_i) &= 1 + X_i \cr
\mu(x_i^{-1}) &= 1 - X_i + X_i^2 - X_i^3 + \cdots \cr
\end{align}

There is a very nice proof that $\mu$ is injective in \cite{KMS}, as well as
discussion of some if its properties.  One such property is that commutators
have zero linear terms.  For example (dropping the $\mu$)

\begin{align}
[x_1,x_2] &= x_1x_2x_1^{-1}x_2^{-1} \cr  &= (1 + X_1)(1 + X_2)(1 - 
X_1 + X_1^2 - \cdots)(1 - X_2 + X_2^2 - \cdots) \cr &= 1 + X_1X_2 - X_2X_1 +
O(3).
\end{align}

Now there is a fairly obvious ordering of 
$Z[[X_1, \dots,X_n]]$.  Write a power series in ascending degree, and 
within each degree list the monomials lexicographically according to subscripts). 
Given two series, order them according to the coefficient of the first term (when
written in the standard form just described) at which they differ.  Thus, for example, $1$ and
$[x_1,x_2]$ first differ at the $X_1X_2$ term, and we see that $1 < [x_1,x_2]$.  It
is not difficult to verify that this ordering, restricted to the group $\{ 1
+O(1)\}$ is invariant under both left- and right-multiplication. 

\begin{exer}
Write $x_1$, $x_2x_1x_2^{-1}$ and $x_2^{-1}x_1x_2$ in increasing order, according to the ordering just described.
\end{exer}

Now we define an ordering of $P_n$.  If $\a$ and $\b$ are pure $n$-braids, compare their Artin coordinates (as described earlier)
$$(\a_1,\a_2,\dots,\a_{n-1}) {\rm\ and\ } (\b_1,\b_2,\dots,\b_{n-1})$$
lexicographically, using within each $F_k$ the Magnus ordering described above.
This all needs choices of conventions, for example, for generators of the free
groups, described in detail in \cite{KR}.  The crucial fact is that the action associated with the semidirect product, by automorphisms $\varphi$, has the property mentioned in Lemma \ref{orderfree}.  We recall the definition of a
positive braid according to Garside: a braid is Garside-positive if it can be
expressed as a word in the standard generators $\s_i$ with only positive
exponents.

\begin{thm}
{\bf (Kim-Rolfsen)} $P_n$ has a bi-ordering with the property that
Garside-positive pure braids are greater than the identity, and the set of all
Garside-positive pure braids is well-ordered by the ordering.
\end{thm}

{\bf Algebraic consequences:} The orderability of the braid groups has
implications beyond what we already knew -- e. g. that they are torsion-free.  In
the theory of representations of a group $G$, it is important to understand the
group algebra $\C G$ and the group ring $\Z G$.  These rings also play a role in
the theory of Vassiliev invariants.  A basic property of a ring would be whether
it has (nontrivial) zero divisors.   If a group $G$ has an element $g$ of finite order, say 
$g^p = 1$, but no smaller power of $g$ is the identity, then we can calculate in $\Z G$
$$(1 - g)(1 + g + g^2 + \cdots + g^{p-1}) = 1 - g^p = 0$$
and we see that both terms of the left-hand side of the equation are (nonzero) divisors
of zero.  A long-standing question of algebra is whether the group ring $\Z G$ can have 
zero divisors if $G$ is torsion-free.  We do know the answer for orderable groups: 

\begin{exer}
If $R$ is a ring without zero divisors, and $G$ is a
right-orderable group, then the group ring $RG$ has no zero divisors.  Moreover,
the only units of $RG$ are the monomials $rg$, $g \in G, r$ a unit of $R$.
\end{exer}

\begin{prop}
{\bf (Malcev, Neumann)} If $G$ has a bi-invariant ordering, then its group ring $\Z G$ embeds in a division algebra, that is, an
extension in which all nonzero elements have inverses.
\end{prop}

These results give us new information about the group rings of the braid
groups.

\begin{thm}  
$\Z B_n$ has no zero divisors.  Moreover,
$\Z P_n$ embeds in a division algebra. 
\end{thm}

A proof of the theorem of Malcev and Neumann \cite{malcev} 
can be found in \cite{MR}.  

\begin{exer}
Which subgroups of $Aut(F_n)$ are right-orderable?
\end{exer}

Of course, $Aut(F_n)$ itself is not
right-orderable, because it has elements of finite order, e.g. permuting the
generators.

A final note regarding orderings:  As we've seen, the methods we've used for
ordering $B_n$ and $P_n$ are quite different.  One might hope there could be
compatible orderings: a bi-ordering of $P_n$ which extends to a right-invariant
ordering of $B_n$.  But, in recent work with Akbar Rhemtulla \cite{RR}, we showed this is
hopeless!

\begin{thm}
{\bf (Rhemtulla, Rolfsen)}. For $n \ge 5$, there is no right-invariant
ordering of $B_n$, which, upon restriction to $P_n$, is also left-invariant.
\end{thm}

\end{document}